\documentclass{article}
\usepackage{tikz}

\usepackage{mathrsfs}[12.6pt]
\usepackage{bm, amscd,  amsfonts, amssymb, cite,texdraw}
\usepackage{amsmath, epsfig, cite, color, colordvi}
\usepackage{latexsym,graphicx}
\usepackage{caption}
\usepackage{subcaption}

\usetikzlibrary{trees}

\tikzstyle{every node}=[circle,inner sep=1pt,fill=white!60]
\tikzstyle{tn}=[shape=circle, draw, color=black!70]
\tikzstyle{marke}=[shape=circle,minimum size=0.2cm, draw,blue]

\makeatletter \@addtoreset{equation}{section} \makeatother

\newtheorem{theo}{Theorem}[section]
\newtheorem{thm}[theo]{Theorem}

\newtheorem{lem}[theo]{Lemma}

\def\0{\mathbb 0}

\def\gen{{\rm Gen}}
\def\deg{{\rm deg}}

\def\pf{\noindent {\it Proof.\ }}
\def\qed{\hfill \rule{4pt}{7pt}}
\begin{document}

\begin{center}
{\Large\bf  A Context-free Grammar for  \\[10pt]
 the Ramanujan-Shor Polynomials }

\vskip 6mm

 William Y.C. Chen$^1$ and Harold R.L. Yang$^2$\\[6pt]
 $^1$Center for Applied Mathematics\\
 Tianjin University\\
 Tianjin 300072, P. R. China\\[9pt]
   $^2$School of Science\\
Tianjin University of Technology and Education\\
Tianjin 300222, P. R. China

\vskip 3mm

 E-mail: $^1$chenyc@tju.edu.cn, $^2$yangruilong@mail.nankai.edu.cn

\end{center}

\vskip 6pt

\centerline{Dedicated to Joseph P. S. Kung}

\vskip 3pt

\begin{abstract}
 Ramanujan defined the polynomials $\psi_{k}(r,x)$ in his study of
 power series inversion. Berndt, Evans and Wilson obtained a recurrence
 relation for $\psi_{k}(r,x)$. In a different context,
Shor introduced the polynomials $Q(i,j,k)$ related to
 improper edges of a rooted tree, leading to a refinement of
 Cayley's formula. He also proved a recurrence relation
 and raised the question of finding a combinatorial proof.
  Zeng realized that
 the polynomials of Ramanujan coincide with
 the polynomials of Shor, and that the recurrence relation of Shor coincides with the recurrence relation of Berndt, Evans and Wilson.
  So we call these polynomials the Ramanujan-Shor polynomials, and call the  recurrence relation the Berndt-Evans-Wilson-Shor recursion.
 A combinatorial proof of this recursion was obtained by Chen and Guo, and a simpler proof was recently given by Guo.
From another perspective,
Dumont and Ramamonjisoa found a context-free grammar $G$
to generate the number of rooted trees on $n$ vertices with $k$
improper edges.
Based on the grammar $G$, we find a
grammar $H$  for the Ramanujan-Shor polynomials. This leads to
a formal calculus for the Ramanujan-Shor polynomials.  In particular,
 we obtain a grammatical derivation of the  Berndt-Evans-Wilson-Shor recursion.
We also provide a grammatical approach to the Abel identities
and a grammatical explanation of the Lacasse identity.
\end{abstract}



\begin{section}{Introduction}

For integers $1\leq k\leq r+1$, Ramanujan \cite{Ramanujan} defined the polynomials $\psi_k(r,x)$ by the following relation:
\begin{equation}\label{equation:psi}
  \sum_{k=0}^\infty \frac{(x+k)^{r+k}e^{-u(x+k)}u^{k}}{k!}
  =\sum_{k=1}^{r+1}\frac{\psi_k(r,x)}{(1-u)^{r+k}},
\end{equation}
and derived the recurrence relation:
\begin{equation}\label{equation:relation_original}
  \psi_k(r,x)=(x-1)\psi_k(r,x-1)+\psi_{k-1}(r+1,x)-\psi_{k-1}(r+1,x-1),
\end{equation}
where  $\psi_1(0,x)=1$, $\psi_0(r,x)=0$ and
$ \psi_k(r,x)=0$ for $k>r+1$.
Berndt, Evans and Wilson \cite{Berndt} obtained
another  recurrence relation
\begin{equation} \label{Berndt}
\psi_k(r,n)=(n-r-k+1)\psi_k(r-1,n)+(r+k-2)\psi_{k-1}(r-1,n).
\end{equation}
Ramanujan also proved the identity
\begin{equation}\label{equation:psi_sum}
\sum_{k=1}^{r+1}\psi_k(r,x)=x^r.
\end{equation}

Zeng \cite{Zeng} observed that the polynomials
 $\psi_k(r,x)$  coincide with the polynomails
 introduced by Shor \cite{Shor}.
Let
\begin{equation} \label{Qnkx}
Q_{n,k}(x)=\psi_{k+1}(n-1,x+n).
\end{equation}
Then \eqref{equation:psi_sum} can be rewritten as
\begin{equation}\label{Qnkx-sum}
\sum_{k=0}^{n-1}Q_{n,k}(x)=(x+n)^{n-1},
\end{equation}
and \eqref{Berndt} can be recast as
\begin{equation}\label{equation:relation_question}
   Q_{n,k}(x)=(x-k+1)Q_{n-1,k}(x+1)+(n+k-2)Q_{n-1,k-1}(x+1),
\end{equation}

 Shor \cite{Shor} defined the polynomials $Q_{n,k}(x)$ in the notation  $Q(i,j,k)$.
 To be more specific, the number $Q(i,j,k)$ is defined by the recurrence relation
 \begin{equation}\label{Q_rec}
Q(i,j,k)=(k+i-1)Q(i-1,j,k)+(i+j-2)Q(i-1,j-1,k),
\end{equation}
where $Q(1,0,k)=1$, $Q(i,-1,k)=0$ for $i\geq1$ and $Q(1,j,k)=0$ for $j\geq1$.

Let  $[n]$ denote the set $\{1,2,\ldots,n\}$. Shor \cite{Shor} showed that,
for $k\geq 1$,   $kQ(i-k,j,k)$ equals the number of forests on $[i]$ rooted
at $\{1,2,\ldots,k\}$ with $j$ improper edges, where an improper edge of a rooted tree $T$ is an edge $(u,v)$ in $T$ such that there exists a descendant of $v$ that is smaller than $u$.
This leads to a combinatorial proof of the relation
\begin{equation}\label{equation:sum_Shor}
  \sum_{j=0}^{i-1}Q(i,j,k)=(i+k)^{i-1},
\end{equation}
 which can be considered as a refinement of Cayley's formula.
 As noted by Shor \cite{Shor}, $Q(i,j,k)$ is a polynomial in $k$ for $i\geq1$ and $0\leq j\leq i-1$, and \eqref{equation:sum_Shor} holds when $k$ is replaced by a variable $x$. Indeed, the polynomials $Q(i,j,x)$ can be expressed as
 $Q_{i,j}(x)$, and so we call $Q_{n,k}(x)$ the Ramanujan-Shor polynomials.

 Shor proved that in addition to the recurrence relation \eqref{Q_rec},
 $Q(i,j,k)$ also satisfies the following relation: For $i\geq1$, $0\leq j\leq i-1$ and $k\geq1$,
 \begin{equation}\label{eq:equ}
Q(i,j,k)=(k-j+1)Q(i-1,j,k+1)+(i+j-2)Q(i-1,j-1,k+1).
\end{equation}
 He asked the question of
 finding a combinatorial interpretation of the above recurrence relation. Notice that
 \eqref{eq:equ} is precisely the recursion \eqref{equation:relation_question} proved by Berndt, Evans and Wilson. So we call it the   Berndt-Evans-Wilson-Shor recursion.

It is worth mentioning that the first recursion
\eqref{Q_rec} of Shor can be deduced from the recursion \eqref{equation:relation_original} of Ramanujan and the
Berndt-Evans-Wilson-Shor recursion \eqref{equation:relation_question}.

Zeng \cite{Zeng} found the following interpretations of the polynomials $Q_{n,k}(x)$ in terms of the
 number of improper edges of trees on $[n+1]$ with root $1$:
\begin{equation}\label{zeng_1}
Q_{n,k}(x)=\sum_{T\in F_{n+1,k}}x^{\deg_T(1)-1},
\end{equation}
where $F_{n,k}$ denotes the set of  trees on $[n]$ with $k$ improper edges and with root $1$, and  $\deg_T(1)$ denotes the degree of the vertex $1$ in $T$.  Zeng also showed that the polynomials $Q_{n,k}(x)$ can   be
interpreted by the number of improper edges of rooted trees (not necessarily rooted at $1$) on $[n]$, namely,
\begin{equation}\label{zeng_2}
Q_{n,k}(x)=\sum_{T\in R_{n,k}}(x+1)^{\deg_T(1)},
\end{equation}
where $R_{n,k}$ denotes the set of rooted trees on $[n]$ with $k$ improper edges. In answer to a question of Zeng \cite{Zeng},
 Chen and Guo \cite{CG} found a  bijection
   showing that \eqref{zeng_1} and \eqref{zeng_2} are equivalent.

We should also note that \eqref{zeng_1} is equivalent to the
 interpretation of $Q_{n,k}(x)$ given by Shor when  $x$ is a positive integer.  As noted by Shor \cite{Shor}, for a positive integer $r$,  $rQ_{n,k}(r)$ equals the number of forests on $[n+r]$ rooted
at $\{1, 2, \ldots, r\}$ with a total number of $k$ improper edges.
Let $F$ be such a forest counted by $rQ_{n,k}(r)$.
Let $T_i$ be the tree in $F$ rooted at $i$, where $1\leq i \leq r$.
For each $T_i$, removing the root $i$ and coloring the subtrees of $T_i$ with color $i$, we get a forest on $\{r+1, r+2, \ldots, r+n\}$ with each tree
colored by one of  colors $1, 2, \ldots, r$. After relabeling, this
leads to a forest on $[n]$ with each tree associated with one
of the colors $1, 2, \ldots, r$.
Let $U_{n,k}$ denote the set of forests of rooted trees on $[n]$ with $k$
 improper edges. For a forest $F$ in $U_{n,k}$, let ${\rm tree}(F)$ denote the number of trees in $F$. By the above argument, one sees that
\begin{equation*}
    rQ_{n,k}(r)=\sum_{F\in U_{n,k}}r^{{\rm tree}(F)},
\end{equation*}
which is equivalent to \eqref{zeng_1} since a forest $F$ in $U_{n,k}$
 gives rise to a rooted tree $T$ in $F_{n+1,k}$ by adding a new root $0$.


Subtracting \eqref{equation:relation_question} from
 \eqref{Q_rec},  the Berndt-Evans-Wilson-Shor recursion \eqref{equation:relation_question} takes the form
\begin{equation}\label{eq:rq_equ}
Q_{n,k}(1+x)=Q_{n,k}(x)+(n+k-1)Q_{n-1,k}(1+x),
\end{equation}
where $n\geq 1, 0\leq k \leq n-1$,   $Q_{1,0}(x)=1$ and
$Q_{n,k}(x)=0$ if $k\geq n$ or $k<0$.
 Chen and Guo \cite{CG}  gave a combinatorial proof of \eqref{eq:rq_equ} in answer to the question of Shor. More precisely, let  $T_{n+1,k}[\deg(2)>0]$ denote the set of rooted trees in $T_{n+1,k}$ for which the vertex $2$ is not a leaf, and let $T_{n+1,k}[\deg(n+1)>0]$  denotes  the set of rooted trees in $T_{n+1,k}$ for which the vertex $n+1$ is not a leaf.
 A bijection between $T_{n+1,k}[\deg(2)>0]$ and $T_{n+1,k}[\deg(n+1)>0]$ was
 constructed in \cite{CG}.
   A simpler bijection was given by Guo \cite{Guo}.

Based on Shor's recursive procedure to construct rooted trees, Dumont and Ramamonjisoa \cite{Dumont} found a context-free grammar to enumerate rooted trees with a given number of improper edges. They defined a
grammar $G$ by the following substitution rules:
\[
G\colon A\rightarrow A^3S,\; S\rightarrow AS^2.
\]
Let $D$ denote the formal derivative with respect to $G$.
Dumont and Ramamonjisoa showed that, for $n\geq1$,
\begin{equation*}
   D^{n-1}(AS)=A^nS^{n}\sum_{k=0}^{n-1}b(n,k)S^{k},
\end{equation*}
where $b(n,k)$ denotes the number of rooted trees on $[n]$ with $k$ improper edges. Note that $b(n,k) =Q_{n,k}(0)$.

Based on the Dumont-Ramamonjisoa grammar, we obtain a grammar $H$ to generate the Ramanujan-Shor polynomials $Q_{n,k}(x)$.
Let
\[
H\colon a\rightarrow axy,x\rightarrow xyw,y\rightarrow y^3w,t\rightarrow yw^2,
\]
and let $D$  denote the formal derivative with respect to $H$.
For $n\geq1$, we obtain the following relation
\[
  D^{n}(a)=axy^nw^{n-1}\sum_{k=0}^{n-1}Q_{n,k}(xw^{-1})y^k.
\]
With the aid of the grammar $H$, we are led to
a simple derivation of the Berndt-Evans-Wilson-Shor recursion in the
form of \eqref{eq:rq_equ}.

It turns out that the grammar $H$ can also be used to deal with the
 Abel identities. In  a certain sense, the formal derivative
 with respect to the grammar $H$ can be
 viewed as a shift-invariant operator
 for the Abel identities in the spirit of the umbral calculus, see Rota \cite{Rota}. As will be seen, the Abel identities can be
 deduced from the Leibnitz formula with respect to the grammar $H$.

 Riordan \cite {Riordan} defined the sum
 \[
A_n(x_1,x_2;p,q)=\sum_{k=0}^n\binom{n}{k}(x_1+k)^{k+p}(x_2+n-k)^{n-k+q},
\]
where $n\geq 1$ and the parameters $p,q$ are integers.
He found closed formulas of $A_n(x_1,x_2;p,q)$ for some
$p$ and $q$. These identities were called the
Abel identities or the Abel-type identities since  the case   $(p,q)=(-1,0)$  corresponds to  the classical Abel identity. We give a grammar $H'$ based on the grammar $H$ and show that
the summations $A_n(x_1, x_2; p,q)$ can be evaluated by
using the grammar $H'$.
Using this approach, closed forms can be deduced  for  $A_n(x,y;-1,0)$, $A_n(x_1,x_2,-1,-1)$ and $A_n(x_1,x_2,-2,0)$ and $A_n(x_1,x_2;-2,-2)$. The case for
  $A_n(x_1,x_2;-2,-2)$ seems to be new.

We conclude this paper with a simple grammatical explanation of the
identity:
  \begin{equation*}
  n^{n+1}=
  \sum_{k=1}^n\sum_{k=0}^{n-j}\binom{n}{j}\binom{n-j}{k}j^jk^k(n-j-k)^{n-j-k}.
  \end{equation*}
We call this identity the Lacasse
identity. It was conjectured by Lacasse \cite{Lacasse} in the study of the
PAC-Bayesian machine learning theory. Since then, several proofs have been found. For example, Sun \cite{Sun} gave a derivation by using the umbral calculus, Younsi \cite{Younsi} found a proof with aid of the Abel identity, Prodinger \cite{Prodinger} provided a justification  based on Cauchy's integral formula, Gessel \cite{Gessel} proved the identity by means of the Lagrange inversion formula, and Chen, Peng and Yang \cite{CPY} obtained a combinatorial
interpretation in terms of triply rooted trees.

This paper is organized as follows.
In Section \ref{section:Dumont}, we give an overview of the Dumont-Ramamonjisoa grammar and introduce a grammatical labeling of labeled trees.
In Section \ref{section:general_g}, we find a grammar $H$ to  generate the Ramanujan-Shor polynomials.
 Section \ref{section:relation} is devoted to a proof of the Berndt-Evans-Wilson-Shor relation by using the grammar $H$.
In Section \ref{section:Abel}, we consider the grammatical derivations of  Abel identities. We also provide a  grammatical explanation of the Lacasse identity.

\end{section}

\section{The Dumont-Ramamonjisoa  Grammar }\label{section:Dumont}

In this section, we give an overview of the context-free grammar introduced by Dumont and Ramamonjisoa \cite{Dumont} to generate rooted trees.
The approach of using context-free grammars to
combinatorial polynomials was introduced in \cite{Chen}. Further studies can be found in \cite{Fu,dumont,Dumont,MMY,MY}.  A context-free grammar $G$ over an  alphabet $A$ is defined to be a set of production rules.
 Given a context-free grammar,   one may define a formal derivative $D$ as a
 differential operator  on polynomials or Laurent polynomials in  $A$, that is, $D$ is a linear operator satisfying the relation
  \[
  D(uv)=D(u)v+uD(v),
  \]
  and in general  the Leibnitz formula
  \begin{equation}\label{equation:Leibnitz}
    D^n(uv)=\sum_{k=0}^n\binom{n}{k}D^k(u)D^{n-k}(v).
  \end{equation}

 Dumont and Ramamonjisoa \cite{Dumont} defined the following grammar \begin{equation}\label{Ramanujan_grammar}
G\colon A\rightarrow A^3S,S\rightarrow AS^2.
\end{equation}
Let $D$ denote the formal derivative with respect to the
grammar $G$. Notice that $D$ can also be viewed as the operator
\[
D=A^3S\frac{\partial}{\partial A}+AS^2\frac{\partial}{\partial S}.
\]
Dumont and Ramamonjisoa established a connection between the grammar $G$ and the enumeration of rooted trees on $[n]$ with $k$ improper edges. The notion
of an improper edge of a rooted tree was introduced by Shor.
 Let $T$ a rooted tree on $[n]$. An edge of $T$ is represented
  by a pair $(u,v)$ of vertices with $v$ being a child
  of $u$. We say that an edge $(u,v)$ of $T$ is improper  if there exists a descendant of $v$ that is smaller than $u$.
 Let $b(n,k)$ denote the number of rooted trees on $[n]$ with $k$ improper edges. Dumont and Ramamonjisoa obtained the following relation.

\begin{thm}\label{th:Dumont}
  For $n\geq1$,
  \begin{equation}\label{eq:Dumont1}
    D^{n-1}(AS)=A^nS^{n}\sum_{k=0}^{n-1}b(n,k)A^k.
  \end{equation}
\end{thm}

For example, for $n=1,2,3$, we have
\begin{align*}
  D^0(AS)&=AS,\\[6pt]
  D^1(AS)&=D(A)S+AD(S)=A^2S^2(1+A),\\[6pt]
  D^2(AS)&=D(D(AS))=A^3S^3(2+4A+3A^2).
\end{align*}

Dumont and Ramamonjisoa gave a proof of the above theorem by showing that the coefficients of $D^n(AS)$ satisfy the recurrence relation \eqref{Q_rec} of Shor.
More precisely, let $s(n,k)$ denote the coefficient of $A^{n+k}S^n$ in $D^{n-1}(AS)$, Dumont and  Ramamonjisoa proved that
\[
s(n,k)=(n-1)s(n-1,k)+(n+k-2)s(n-1,k-1),
\]
which is equivalent to the relation \eqref{Q_rec}
for the case  $x=0$.

Here we present a proof in the language of a grammatical
labeling of rooted trees, which was introduced in
\cite{Fu}.  Let $R_n$ denote the set of rooted trees on $[n]$£¬
and let $F_n$ denote the set of rooted trees on $[n]$
 with $1$. For a rooted tree $T\in R_n$ and a vertex $u$ in $T$, let $\beta_T(u)$ or simply $\beta(u)$ denote the minimum vertex among the descendants of $u$ in $T$. An improper edge $e=(u,v)$ of $T$ is defined as an edge such that $u>\beta(v)$. Otherwise, we say that $e$ is a proper edge.

Let $F_{n,k}$ denote the set of rooted trees in $F_n$ with $k$ improper edges, and let $R_{n,k}$ denote the set of rooted trees in $R_n$ with $k$ improper edges. Shor \cite{Shor} presented a construction of a rooted tree in $ R_{n}$
by adding the vertex $n$ into a tree  in $ R_{n-1}$.
For a better understanding of the construction,
let us consider the following
 procedure to delete the vertex $n$ from  a
 rooted tree $T$ in $R_{n}$ to obtain a
 rooted tree $T'$ in $R_{n-1}$:

\begin{enumerate}
\item Case 1: $n$ is a leaf in $T$. Delete the vertex $n$.
\item Case 2: $n$ is not a leaf. Assume that $n$ has $t$ children $b_1,b_2,\ldots,b_t$. We may further assume that
    \[
    \beta(b_1)<\beta(b_2)<\cdots<\beta(b_t).
     \]
     Contract the edge $(n,b_t)$ and relabel the resulting vertex by $b_t$.
\end{enumerate}

Conversely, one can construct a rooted tree $T$ on $[n]$ with $k$ or $k+1$ improper edges from a rooted tree $T'$ on  $[n-1]$ with $k$ improper edges. There are four operations to add
the vertex $n$ to $T'$.

\begin{enumerate}
\item Adding $n$ to the tree $T'$ as a child of an arbitrary vertex $v$, we obtain a tree $T\in R_{n,k}$ with $n$ being a leaf.
\item Splitting a proper edge $(i,j)$ into $(i,n)$ and $(n,j)$, we obtain a tree $T\in R_{n,k+1}$. In this case, the degree of $n$ equals one.
\item Splitting an improper edge $(i,j)$ into $(i,n)$ and $(n,j)$, we also obtain a tree $T\in R_{n,k+1}$. In this case, the degree of $n$ also equals one.
\item Choose an improper edge $(v,b_j)$ in $T'$, where $v$ has $t$ children $b_1,b_2,\ldots,b_t$ listed in the increasing order of their
    $\beta$-values.
    We relabel $v$ by $n$ and make $v$ a child of $n$. Moreover, assign  $b_1,\ldots,b_j$ to be the children of $n$ and assign $b_{j+1},\ldots,b_t$ to be the children of $v$. Then we are led to a tree $T\in R_{n,k+1}$. In this case, the degree of $n$ in $T$ is
     at least two.
\end{enumerate}

As will be seen, the above construction is
closely related to the grammar $G$. To
demonstrate this connection,
we introduce a grammatical labeling of rooted trees.
We may view a rooted tree $T$ on $[n]$ as a rooted tree $\hat{T}$ on
$\{0, 1, 2, \ldots, n\}$ with $0$ being the root with only one child.
Clearly, the edge below the root $0$ of $\hat{T}$
is always a proper edge.
Moreover, we represent an improper edge by double edges,  called the left edge and the right edge.
 The idea of using double edges to represent an improper
edge is due to Dumont and Ramamonjisoa \cite{Dumont}.
We label a vertex of $\hat{T}$ except for $0$ by $S$ and label an edge of $\hat{T}$ by $A$.
In other words, 
 a proper edge of ${T}$ is labeled by $A$ and an improper edge of ${T}$ is labeled by $A^2$.
The weight of $T$ is defined by the product of the labels attached to $\hat{T}$,
denoted by $w(T)$.
Apparently, for any tree $T$ in $R_{n,k}$, we have $w(T)=A^{n+k}S^n$.

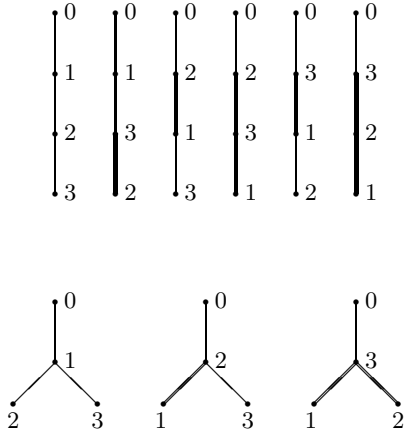
\begin{figure}[!ht]
\begin{center}
\begin{picture}(125,20)
\setlength{\unitlength}{0.8mm}

\multiput(0,15)(10,0){6}{\circle*{1}}
\multiput(0,-15)(10,0){6}{\circle*{1}}
\multiput(0,-5)(10,0){6}{\circle*{1}}
\multiput(0,5)(10,0){6}{\circle*{1}}
\multiput(0,5)(10,0){6}{\line(0,1){10}}
\multiput(0,-15)(20,0){3}{\line(0,1){10}}
\multiput(9.8,-15)(20,0){3}{\line(0,1){10}}\multiput(10.2,-15)(20,0){3}{\line(0,1){10}}
\multiput(0,-5)(10,0){2}{\line(0,1){10}}
\multiput(19.8,-5)(10,0){4}{\line(0,1){10}}\multiput(20.2,-5)(10,0){4}{\line(0,1){10}}
\multiput(1.5,14)(10,0){6}{\small{$0$}}
\multiput(1.5,4)(10,0){2}{\small{$1$}}\multiput(21.5,4)(10,0){2}{\small{$2$}}
\multiput(41.5,4)(10,0){2}{\small{$3$}}\multiput(1.5,-6)(50,0){2}{\small{$2$}}
\multiput(11.5,-6)(20,0){2}{\small{$3$}}
\multiput(21.5,-6)(20,0){2}{\small{$1$}}
\multiput(1.5,-16)(20,0){2}{\small{$3$}}\multiput(11.5,-16)(30,0){2}{\small{$2$}}
\multiput(31.5,-16)(20,0){2}{\small{$1$}}

\multiput(1.5,-34)(25,0){3}{\small{$0$}}
\multiput(0,-33)(25,0){3}{\circle*{1}}
\multiput(0,-43)(25,0){3}{\circle*{1}}\multiput(0,-43)(25,0){3}{\line(0,10){10}}
\multiput(-7,-50)(25,0){3}{\circle*{1}}\multiput(7,-50)(25,0){3}{\circle*{1}}
\put(-7,-50){\line(1,1){7}}
\multiput(17.8,-50)(25,0){2}{\line(1,1){7}}\multiput(18.2,-50)(25,0){2}{\line(1,1){7}}
\multiput(7,-50)(25,0){2}{\line(-1,1){7}}
\put(56.8,-50){\line(-1,1){7}}\put(57.2,-50){\line(-1,1){7}}
\put(1.5,-44){\small{$1$}}\put(26.5,-44){\small{$2$}}\put(51.5,-44){\small{$3$}}
\multiput(-8,-54)(64,0){2}{\small{$2$}}\multiput(6,-54)(25,0){2}{\small{$3$}}
\multiput(16.5,-54)(25,0){2}{\small{$1$}}
\end{picture}
\vspace{4.5cm}
\end{center}
\caption{Rooted trees in $R_3$.}
\label{figure_label}
\end{figure}

Figure \ref{figure_label} illustrates all rooted trees
on $\{1,2, 3\}$, where the improper edges are represented by
double edges, and the vertex $0$ is added at the top of each tree in
$R_3$.
\vskip 1cm

The following relation is a restatement of Theorem \ref{th:Dumont}.

\begin{thm}\label{grammar_w:Dnxy}
  For $n\geq1$,
  \begin{equation}\label{equation:dumont}
  D^{n-1}(AS)=\sum_{T\in R_n}w(T).
  \end{equation}
\end{thm}

In view of the above grammatical laleling of rooted trees, it can be seen that
 the four cases in Shor's construction of a tree $T'$ on $[n]$ from
 a tree on $[n-1]$ correspond to the substitution rules in $G$.
 Instead of giving a detail proof, let us use an example to demonstrate
  the correspondence.

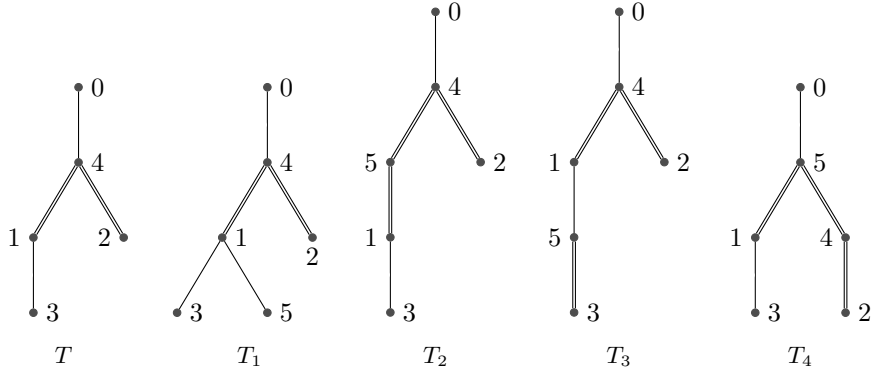
\begin{figure}[!ht]
\subcaptionbox*{$T$}[.195\textwidth]{%
\begin{tikzpicture}
\node [tn,label=0:$0$]{}[grow=down]
	[sibling distance=12mm,level distance=10mm]
    child {node [tn,label=0:{$4$}](four){}
        child {node [tn,label=180:{$1$}](one){}
	       child{node [tn,label=0:{$3$}]{}
                edge from parent
                    node[left ]{}
            }}
        child {node [tn,label=180:{$2$}](two){}}
        edge from parent
            node[left ]{}
        };
\draw [double](four)--node[left=2pt]{}(one);
\draw[double](four)--node[right=2pt]{}(two);
\end{tikzpicture}
}
\subcaptionbox*{$T_1$}[.195\textwidth]{%
\begin{tikzpicture}
\node [tn,label=0:$0$]{}[grow=down]
	[sibling distance=12mm,level distance=10mm]
    child {node [tn,label=0:{$4$}](four){}
        child {node [tn,label=0:{$1$}](one){}
	       child{node [tn,label=0:{$3$}]{}
                edge from parent
                    node[left=1pt ]{}
            }
           child{node [tn,label=0:{$5$}]{}
                edge from parent
                    node[right=1pt ]{}
            }
        }
        child {node [tn,label=270:{$2$}](two){}}
        edge from parent
            node[left ]{}
        };
\draw [double](four)--node[left=2pt]{}(one);
\draw[double](four)--node[right=2pt]{}(two);
\end{tikzpicture}
}
\subcaptionbox*{$T_2$}[.195\textwidth]{%
\begin{tikzpicture}
\node [tn,label=0:$0$]{}[grow=down]
	[sibling distance=12mm,level distance=10mm]
    child {node [tn,label=0:{$4$}](four){}
        child {node [tn,label=180:{$5$}](five){}
            child {node [tn,label=180:{$1$}](one){}
	       child{node [tn,label=0:{$3$}]{}
                edge from parent
                    node[left ]{}
            }}
        }
        child {node [tn,label=0:{$2$}](two){}}
        edge from parent
            node[left ]{}
        };
\draw [double](four)--node[left=2pt]{}(five);
\draw [double](one)--node[left=2pt]{}(five);
\draw[double](four)--node[right=2pt]{}(two);
\end{tikzpicture}
}
\subcaptionbox*{$T_3$}[.19\textwidth]{%
\begin{tikzpicture}
\node [tn,label=0:$0$]{}[grow=down]
	[sibling distance=12mm,level distance=10mm]
    child {node [tn,label=0:{$4$}](four){}
        child {node [tn,label=180:{$1$}](five){}
            child {node [tn,label=180:{$5$}](one){}
	       child{node [tn,label=0:{$3$}](three){}
                edge from parent
                    node[left ]{}
            }}
        }
        child {node [tn,label=0:{$2$}](two){}}
        edge from parent
            node[left ]{}
        };
\draw [double](four)--node[left=2pt]{}(five);
\draw [double](one)--node[left=2pt]{}(three);
\draw[double](four)--node[right=2pt]{}(two);
\end{tikzpicture}
}
\subcaptionbox*{$T_4$}[.19\textwidth]{%
\begin{tikzpicture}
\node [tn,label=0:$0$]{}[grow=down]
	[sibling distance=12mm,level distance=10mm]
    child {node [tn,label=0:{$5$}](four){}
        child {node [tn,label=180:{$1$}](one){}
	       child{node [tn,label=0:{$3$}]{}
                edge from parent
                    node[left ]{}
            }}
        child {node [tn,label=180:{$4$}](two){}
            child {node [tn,label=0:{$2$}](five){}}
        }
        edge from parent
            node[left ]{}
        };
\draw [double](four)--node[left=2pt]{}(one);
\draw[double](four)--node[right=2pt]{}(two);
\draw[double](five)--node[right=2pt]{}(two);
\end{tikzpicture}
}
\caption{An example for the operator $D$}
\label{BeginT}
\end{figure}

In
 Figure \ref{BeginT},  $T$ is a rooted tree on $\{1,2,3,4\}$. The  weight
  of $T$ is $w(T)=A^6S^4$. The trees $T_1$, $T_2$, $T_3$ and $T_4$ are obtained from  $T$ in the four cases of Shor's construction.

{Case} $1$: $T_1$ is obtained from $T$ by adding the vertex $5$ as a leaf.
Comparing the weight of $T_1$ and the weight of $T$, it can be seen that
this operation corresponds to the substitution rule $S\rightarrow AS^2$.
Notice that the label $S$ indicates where one can apply this operation.

{Case} $2$: $T_2$ is obtained from $T$  by splitting  the left edge $(4,1)$ into two edges $(4,5)$ and $(5,1)$.
This operation corresponds to the substitution rule $A\rightarrow A^3S$.

{Case} $3$: $T_3$ is obtained from $T$   by splitting  the proper edge $(1,3)$ into  $(1,5)$ and $(5,3)$.
This operation also corresponds to the substitution rule $A\rightarrow A^3S$.

{Case} $4$: $T_4$ is obtained by adding $5$ to $T$ via the following procedure: $4$ is relabeled by $5$, a new vertex $4$ is added as a child of $5$,  the subtree rooted by  $1$ and the subtree rooted by $2$ are assigned as a child of $5$ and a child of  $4$, respectively. It can be seen that this operation also corresponds to the substitution rule $A\rightarrow A^3S$.

The above argument is sufficient to lead to a rigorous proof of  relation \eqref{equation:dumont}.

  \section{A grammar for the R-S polynomials}\label{section:general_g}

In this section, we
 give a grammar $H$ to generate the Ramanujan-Shor polynomials $Q_{n,k}(x)$.
Define
\begin{equation}\label{grammar:gR}
  H: a\rightarrow axy,x\rightarrow xyw,y\rightarrow y^3w,w\rightarrow yw^2.
\end{equation}

Recall that   a rooted tree $T\in F_{n,k}$
is rooted at $1$ and has $k$ improper edges.
For  $T\in F_{n,k}$, we label a proper edge by $y$, and represent each improper edge of $T$ by double edges, each of which is
 labeled by $y$.
Meanwhile, we label the vertex $1$ by  $a$, label each  child of the vertex $1$ by
 $x$ and label other vertices by $w$, so that for $T\in F_{n+1,k}$,
 the weight of $T$ equals
\begin{equation}\label{weight_rep}
w(T)=ax^{\deg_T(1)}y^{n+k}w^{n-\deg_T(1)}.
\end{equation}

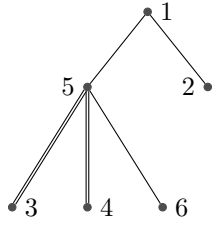
\begin{figure}[!ht]
\begin{center}
\begin{tikzpicture}
\node [tn,label=0:{$1$}]{}[grow=down]
	[sibling distance=16mm,level distance=10mm]
        child {node [tn,label=180:{$5$}](five){}
            [sibling distance=10mm,level distance=16mm]
	       child{node [tn,label=0:{$3$}](three){}}
            child{node [tn,label=0:{$4$}](four){}}
            child{node [tn,label=0:{$6$}]{}
                edge from parent
                    node[right=2pt ]{}
            }
            edge from parent
                node[left=2pt ]{}
            }
        child {node [tn,label=180:{$2$}]{}
            edge from parent
                node[left=2pt ]{}
                }
        ;
\draw [double](four)--(five);
\draw[double](five)--(three);
\end{tikzpicture}
\end{center}
\caption{A rooted tree $T\in F_{6,2}$}\label{Figure:Example}
\end{figure}

Figure \ref{Figure:Example} illustrates a rooted tree in $F_{6,2}$ with weight  $w(T)=ax^2y^7w^3$.

Let $D$ denote the formal derivative with respect to the grammar $H$. Recall that $F_n$ is the set of rooted tree on $[n]$
 with root $1$.
The next theorem shows that $D$ can be used to generate the sum of weights of
rooted trees in $F_n$.

\begin{thm}\label{th:Dna}
For $n\geq1$,
\begin{equation}\label{Ram_ex_Dn}
{D}^n(a)=\sum_{T\in F_{n}}w(T).
\end{equation}
\end{thm}

 To prove the above relation, it is sufficient  to observe that the substitution rules in $H$
 correspond to the changes of  labels in Shor's construction
according to the above labeling scheme.

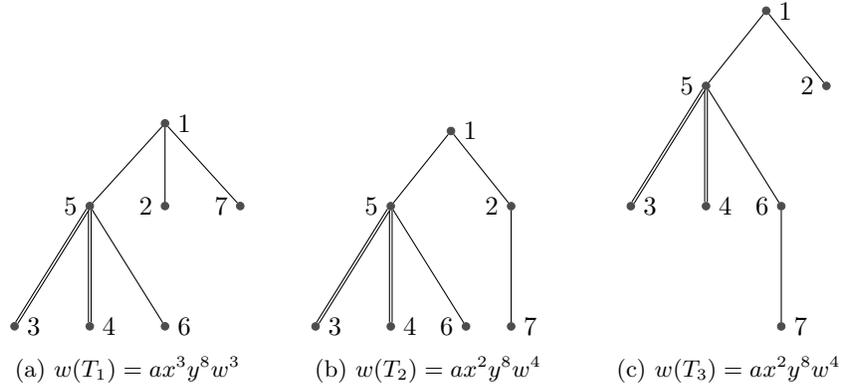
\begin{figure}[!ht]
\begin{center}
\subcaptionbox{$w(T_1)=ax^3y^8w^3$}[.32\textwidth]{%
\begin{tikzpicture}
\node [tn,label=0:{$1$}]{}[grow=down]
	[sibling distance=10mm,level distance=11mm]
        child {node [tn,label=180:{$5$}](five){}
            [sibling distance=10mm,level distance=16mm]
	       child{node [tn,label=0:{$3$}](three){}}
            child{node [tn,label=0:{$4$}](four){}}
            child{node [tn,label=0:{$6$}]{}
                edge from parent
                    node[right=2pt ]{}
            }
            edge from parent
                node[left=2pt ]{}
            }
        child {node [tn,label=180:{$2$}]{}
            edge from parent
                node[left=1pt ]{}
                }
        child {node [tn,label=180:{$7$}]{}
            edge from parent
                node[left=2pt ]{}
                }
        ;
\draw [double](four)--node[left=0.6pt]{}(five);
\draw[double](five)--node[left=3pt]{}(three);
\end{tikzpicture}
}
\subcaptionbox{$w(T_2)=ax^2y^8w^4$}[.32\textwidth]{%
\begin{tikzpicture}
\node [tn,label=0:{$1$}]{}[grow=down]
	[sibling distance=16mm,level distance=10mm]
        child {node [tn,label=180:{$5$}](five){}
            [sibling distance=10mm,level distance=16mm]
	       child{node [tn,label=0:{$3$}](three){}}
            child{node [tn,label=0:{$4$}](four){}
            }
            child{node [tn,label=180:{$6$}]{}
                edge from parent
                    node[left=0.5pt ]{}
            }
            edge from parent
                node[left=2pt ]{}
            }
        child {node [tn,label=180:{$2$}]{}
             [sibling distance=10mm,level distance=16mm]
            child{node [tn,label=0:{$7$}]{}
                edge from parent
                    node [left=2pt]{}
                }
            edge from parent
                node[left=2pt ]{}
                }
        ;
\draw [double](four)--(five);
\draw[double](five)--(three);
\end{tikzpicture}
}
\subcaptionbox{$w(T_3)=ax^2y^8w^4$}[.32\textwidth]{%
\begin{tikzpicture}
\node [tn,label=0:{$1$}]{}[grow=down]
	[sibling distance=16mm,level distance=10mm]
        child {node [tn,label=180:{$5$}](five){}
            [sibling distance=10mm,level distance=16mm]
	       child{node [tn,label=0:{$3$}](three){}}
            child{node [tn,label=0:{$4$}](four){}
            }
            child{node [tn,label=180:{$6$}]{}
                [sibling distance=10mm,level distance=16mm]
                child{node [tn,label=0:{$7$}]{}
                    }
                edge from parent
                    node[left=0.5pt ]{}
            }
            edge from parent
                node[left=2pt ]{}
            }
        child {node [tn,label=180:{$2$}]{}
            edge from parent
                node[left=2pt ]{}
                }
        ;
\draw [double](four)--(five);
\draw[double](five)--node[left=3pt]{}(three);
\end{tikzpicture}
}
\end{center}
\caption{The action of $D$}\label{Figure:Example_D}
\end{figure}

Figure \ref{Figure:Example_D} gives three rooted trees
 $T_1$, $T_2$ and $T_3$ obtained from the tree $T$ in Figure \ref{Figure:Example} by adding the vertex $7$ as a leaf as in Case 1 of Shor's construction. Since $7$ is child of the root $1$,
  $w(T_1)$ is obtained from $w(T)$ by
   applying   the substitution rule $a\rightarrow axy$.
   Similarly, $7$ is a child of $2$ in $T_2$, and $w(T_2)$ is obtained from $w(T)$
   by utilizing the rule $x\rightarrow xyw$.  Since $7$ is a child of $6$ in $T_3$,   $w(T_3)$ is obtained from $w(T)$ by the rule $w\rightarrow yw^2$.

  For Case 2, Case 3 and Case 4 in Shor's construction,
  the changes of weights of the resulting trees can be
  characterized by the rule  $y\rightarrow y^3w$, just like
  the rule $A\rightarrow A^3 S$ in
  the Dumont-Ramamonjisoa grammar.

 We now come to a relationship between the grammar $H$ and  the polynomials  $Q_{n}(x,y)$. Recall that
\[
Q_n(x,y)=\sum_{k=0}^{n-1}Q_{n,k}(x)y^k=\sum_{T\in F_{n+1}}x^{\deg_T(1)-1}y^{\rm imp(T)}.
\]
For $n=1,2,3$, we have
\begin{align*}
    Q_1(x,y)&=1.\\[6pt]
    Q_2(x,y)&=y+x+1.\\[6pt]
    Q_3(x,y)&=3y^2+(3x+4)y+x^2+3x+2.
\end{align*}

\begin{thm}\label{thm:grdna}
For $n\geq1$,
\begin{equation}\label{grammar_cayley_exp}
D^n(a)=axy^nw^{n-1}Q_n(xw^{-1},y).
\end{equation}
\end{thm}

For $n=1,2,3$,  we have
\begin{align*}
  D(a)&=axy=axyQ_0(xw^{-1},y),\\[6pt]
  D^2(a)&=ax^2y^2+axy^2w+axy^3w=axy^2w(xw^{-1}+1+y)=axy^2wQ_1(xw^{-1},y),\\[6pt]
  D^3(a)&=ax^3y^3+3ax^2y^3w+3ax^2y^4w+2axy^3w^2+4axy^4w^2+3axy^5w^2\\[3pt]
  &=axy^3w^2\left(x^2w^{-2}+3xw^{-1}+2+(4xw^{-1}+3)y+3y^2\right)\\[6pt]
  &=axy^3w^2Q_2(xw^{-1},y).
\end{align*}

We end this section with a grammatical derivation of the relation \eqref{Q_rec}.

\begin{thm}
For $n\geq2$ and $1\leq k\leq n-2$, we have
\begin{equation}\label{relation_sim}
Q_{n,k}(x) = (x+n-1)Q_{n-1,k}(x) + (n+ k-2)Q_{n-1,k-1}(x).
\end{equation}
\end{thm}

\pf For $n\geq 1$, by the definition of $Q_n(x,y)$,
\eqref{grammar_cayley_exp} can be written as
\begin{equation}\label{dn_a}
  D^{n}(a)=(xw^{-1})ay^{n}w^{n}\sum_{k=0}^{n-1}y^kQ_{n,k}(xw^{-1}).
\end{equation}
For $n\geq 2$, substituting $n$ by $n-1$,  \eqref{dn_a} takes the form
\begin{equation} \label{dn1a}
  D^{n-1}(a)=(xw^{-1})ay^{n-1}w^{n-1}\sum_{k=0}^{n-2}y^kQ_{n-1,k}(xw^{-1}).
\end{equation}
 Since
 \begin{equation}\label{eq:constantxw}
 D(xw^{-1})=xyw\cdot w^{-1}-x\cdot w^{-2}yw^2=0,
 \end{equation}
 that is, $xw^{-1}$ is a constant with respect to $D$,
 we find that
 \[
 D\left(y^kQ_{n-1,k}(xw^{-1})\right)=Q_{n-1,k}(xw^{-1}) D(y^k) = ky^{k+2}wQ_{n-1,k}(xw^{-1}).
 \]
Meanwhile,
 \[
 D\left(ay^{n-1}w^{n-1}\right)=axy^{n}w^{n-1}+(n-1)ay^{n+1}w^n+(n-1)ay^nw^n.
 \]
Therefore,  applying the operator $D$ to both sides of \eqref{dn1a} yields
 \begin{align}
D^n(a)&=(xw^{-1})ay^nw^n\left\{\sum_{k=0}^{n-2}ky^{k+1}Q_{n-1,k}(xw^{-1})\right.\notag\\[6pt]
&\qquad \left.+(xw^{-1}+(n-1)y+(n-1))\sum_{k=0}^{n-2}y^kQ_{n-1,k}(xw^{-1})\right\}.\label{dn_expand}
\end{align}
For $n\geq2$ and $0\leq k\leq n-2$, comparing the coefficient of $axy^{n+k}w^n$ on the
right hand sides of \eqref{dn_a} and  \eqref{dn_expand},
we deduce that
\begin{equation}
  Q_{n,k}(xw^{-1})=(xw^{-1}+n-1)Q_{n-1,k}(xw^{-1})+(n+k-2)Q_{n-1,k-1}(xw^{-1}).
\end{equation}
Setting $w=1$ completes the proof.
\qed
\section{The Berndt-Evans-Wilson-Shor  Recursion }\label{section:relation}

The section is devoted to a grammatical derivation of the Berndt-Evans-Wilson-Shor  recurrence relation \eqref{equation:relation_question}.
To this end, we establish
a grammatical expression for  $Q_{n}(r+x,y)$, where $r$ is a nonnegative integer.

  \begin{thm}\label{eq:dnaxr_close}
  For $n\geq1$ and $r \geq 0$,
  \begin{equation}\label{eq:dnabr}
  D^n(ax^r)=ax^ry^nw^n(r+xw^{-1})Q_n(r+xw^{-1},y).
  \end{equation}
\end{thm}

To prove this relation, we give a combinatorial interpretation of $Q_n(x+r,y)$
based on Zeng's interpretation of $Q_{n,k}(x)$ in terms of
the set $F_{n,k}$ of rooted trees on $[n]$ with root $1$ containing $k$ improper edges.
We define $F^{(r)}_n$ to be the set of rooted trees on $[n]$ with root 1
for which each child of the root is colored by
one of the colors $b, w_1,w_2,\ldots,w_r$,
where $b$ stands for the black color, and $w_1, w_2, \ldots, w_r$
are considered  white colors.

We now define a grammatical labeling of a rooted tree $\bar{T}\in F^{(r)}_n$.
First,  represent an improper edge of $\bar{T}$ by double edges, and
denote the resulting tree by $\hat{T}$. Then
the root of $\hat{T}$ is labeled by $ax^r$,  a black vertex is
labeled by $x$ and
 each of the remaining vertices is labeled by $w$.
 Moroever, each edge of $\hat{T}$ is labeled by $y$.
In other words, as far as $\bar{T}$ is concerned, a proper edge  is labeled by $y$ and an improper edge is labeled by $y^2$.
For  $\bar{T}\in F^{(r)}_{n}$, we have
\begin{equation}\label{gl:cf}
w(\bar{T})=ax^{{\rm black}(\bar{T})+r}w^{n-{\rm black}(\bar{T})}y^{n+{\rm imp}(\bar{T})},
\end{equation}
where $\rm black(\bar{T})$ denotes the number of black vertices in $\bar{T}$.

Using the above labeling scheme, the right hand side of \eqref{eq:dnabr} can be expressed as follows.

\begin{thm}
  For $n\geq0$ and $r \geq 0$,
  \begin{equation}\label{eq:wFr}
    ax^ry^nw^{n}(r+xw^{-1})Q_n(r+xw^{-1},y)=\sum_{\bar{T}\in F^{(r)}_{n+1}}w(\bar{T}).
  \end{equation}
\end{thm}
\pf By \eqref{gl:cf}, we see that
\begin{align*}
\sum_{\bar{T}\in F^{(r)}_{n+1}}w(\bar{T})&=\sum_{\tilde{T}\in F^{(r)}_{n+1}}ax^{{\rm black}(\bar{T})+r}w^{n-{\rm black}(\bar{T})}y^{n+{\rm imp}(\bar{T})}\\[6pt]
&=ax^ry^nw^{n}\sum_{k=0}^{n-1}y^k\sum_{\bar{T}\in F^{(r)}_{n+1,k}}x^{{\rm black}(\bar{T})}w^{-{\rm black}(\bar{T})}.
\end{align*}
Given a rooted tree $T\in F_{n,k}$, one can construct a rooted tree $\bar{T}$ in $ F^{(r)}_{n,k}$ by assigning the color $b$ to
some children of the root $1$ and  one
of the $r$ white colors to each remaining children of the root $1$. Thus
\begin{align*}
\sum_{\bar{T}\in F^{(r)}_{n+1,k}}x^{{\rm black}(\bar{T})}w^{-{\rm black}(\bar{T})}
&=\sum_{T\in F_{n+1,k}}\sum_{i=0}^{\deg_T(1)}{\deg_T(1)\choose i}(xw^{-1})^{i}r^{\deg_T(1)-i}\\[6pt]
&=\sum_{T\in F_{n+1,k}} (r+xw^{-1})^{\deg_T(1)},
\end{align*}
which can be expressed as
$Q_{n,k}(r+xw^{-1})$ according to the interpretation \eqref{zeng_1}
of $Q_{n,k}(x)$.
It follows that
\begin{equation*}
\sum_{\bar{T}\in F^{(r)}_{n+1}}w(T)= ax^ry^nw^{n}(1+xw^{-1})\sum_{k=0}^{n-1}y^kQ_{n,k}(r+xw^{-1}),
\end{equation*}
as claimed. \qed

The following theorem establishes a connection between
 the grammar $H$ and the sum of weights of rooted trees in $F^{(r)}_n$.

\begin{thm}\label{th:Dnab}
  For $n\geq1$ and $r \geq 0$,
  \[
  D^n(ax^r)=\sum_{\bar{T}\in F^{(r)}_{n+1}}w(\bar{T}).
  \]
  \end{thm}

The proof is similar to that of Theorem \ref{th:Dna}. The operation of adding $n$ as a black child of the root $1$ can be described by the substitution rule $a\rightarrow axy$ and the operation of adding $n$ as a white child of the root $1$ corresponds to the rule $x\rightarrow xyw$.

We now give a grammatical derivation of the Berndt-Evans-Wilson-Shor  recursion for $Q_{n,k}(x)$,  that is, for $n\geq 1$ and $0\leq k \leq n-1$,
  \begin{equation}\label{Q_rec_t}
Q_{n,k}(1+x)=Q_{n,k}(x)+(n+k-1)Q_{n-1,k}(1+x).
\end{equation}
Note that $Q_{1, 0}(x)=1$ and
$Q_{n,k}(x)=0$ if $k\geq n$ or $k<0$.

Our proof relies on the generating function with respect to the grammar $H$. For a Laurent polynomial $w$ of the variables in the alphabet $V$, the exponential generating function of $w$ with respect to $D$ is defined by
\[
\gen(w,t)=\sum_{n\geq0}D^n(w)\frac{t^n}{n!}.
\]
We have the following properties:
\begin{align}
  \gen'(w,t)&=\gen(D(w),t)\\[6pt]
  \gen(w+v,t)&=\gen(w,t)+\gen(v,t)\\[6pt]
  \gen(wv,t)&=\gen(w,t)\gen(v,t),
\end{align}
where $\gen'(w,t)$ stands for the differentiation of
$\gen(w,t)$ with respect to $t$, and $v$ is also a Laurant
polynomial of the variables in the alphabet $V$, see \cite{Chen}.

We are now in a position to present a grammatical proof of \eqref{Q_rec_t}.
It is easily seen that \eqref{Q_rec_t} follows from the following
relation for $n\geq 1$,
\begin{align}
&axy^{n+1}w^n(1+xw^{-1})Q_n(1+xw^{-1},y)\notag\\[6pt]
&\qquad =axy^{n+1}w^n(1+xw^{-1})Q_n(xw^{-1},y)\notag\\[6pt]
&\qquad\qquad +axw^n(1+xw^{-1})\sum_{k=0}^{n-2}(n+k-1)Q_{n-1,k}(1+xw^{-1})y^{n+k+1}.\label{eq:expand}
\end{align}
Invoking \eqref{eq:dnabr} for $n\geq 1$ and $r=0$,  we obtain that for $n\geq 1$,
\begin{equation}\label{eq:dna_expand}
D^n(a)=axy^{n}w^{n-1}Q_n(xw^{-1},y).
\end{equation}
Again, utilizing \eqref{eq:dnabr} for $n\geq 1$ and $r=1$, we find that
\begin{equation}\label{eq:dnax_expand}
D^n(ax)=axy^{n}w^n(1+xw^{-1})Q_n(1+xw^{-1},y),
\end{equation}
and so
\begin{equation}\label{eq:dnax_expand-1}
D^{n-1}(ax)=axy^{n-1}w^{n-1}(1+xw^{-1})Q_{n-1}(1+xw^{-1},y).
\end{equation}
Thus \eqref{eq:expand} can be rewritten as
\begin{equation}\label{eq:geq}
    yD^n(ax)=yw(1+xw^{-1})D^n(a)+y^3w\frac{\partial (D^{n-1}(ax))}{\partial y}.
\end{equation}
Expanding \eqref{eq:dnax_expand-1} as
 \[   D^{n-1}(ax)=axy^{n-1}w^{n-1}(1+xw^{-1})\sum_{k=0}^{n-2}Q_{n-1,k}(1+xw^{-1})y^k,
\]
we see that
\begin{align*}
    axy\frac{\partial (D^{n-1}(ax)) }{\partial a}&=xyD^{n-1}(ax),\\[6pt]
    xyw\frac{\partial (D^{n-1}(ax)) }{\partial x}&=ywD^{n-1}(ax),\\[6pt]
    yw^2\frac{\partial (D^{n-1}(ax)) }{\partial w}&=(n-1)ywD^{n-1}(ax).\\[6pt]
\end{align*}
Notice that
\[
  D= axy\frac{\partial  }{\partial a}+xyw\frac{\partial  }{\partial x}
    +y^3w\frac{\partial  }{\partial y}+yw^2\frac{\partial }{\partial w},
\]
so that
\begin{equation}\label{eq:dnax_new}
    D^n(ax)
    =xyD^{n-1}(ax)+nywD^{n-1}(ax)+y^3w\frac{\partial (D^{n-1}(ax))}{\partial y},
\end{equation}
and therefore, \eqref{eq:geq} is equivalent to
\begin{equation}\label{eq:target_new}
(y-1)D^{n-1}(D(ax))+(nyw+xy)D^{n-1}(ax)=yw(1+xw^{-1})D^{n-1}(D(a)),
\end{equation}
for $n\geq 1$. In terms of the generating functions, \eqref{eq:target_new} can be reformulated as
\begin{align}
(xy+yw)\gen(axy,t)&=(xy+yw)\gen(ax,t)\notag\\[6pt]
&\qquad +(y-1+tyw)\gen(axyw+ax^2y,t).\label{equation:equivalent}
\end{align}
Let
\begin{align*}
  A(t)&=(y-1+tyw)\gen(axyw+ax^2y,t)\\[6pt]
  &\qquad +(xy+yw)\gen(ax,t)-(xy+yw)\gen(axy,t).
\end{align*}
Since $D(xw^{-1})=0$ as given in \eqref{eq:constantxw}, we have
\[
A(t)=(1+xw^{-1})\gen(axyw,t)
\left(y-1+tyw+yw\gen(y^{-1}w^{-1}-w^{-1},t)\right).
\]
It remains to show that
\begin{equation}\label{eq:gen_eq}
y-1+tyw+yw\gen(y^{-1}w^{-1}-w^{-1},t)=0.
\end{equation}
Observe that
\begin{equation*}
D(y^{-1}w^{-1}-w^{-1})=-y^{-2}w^{-1}y^3w-y^{-1}w^{-2}yw^2+w^{-2}yw^2=-1.
\end{equation*}
Hence
\begin{equation}\label{equation:x-1}
\gen(y^{-1}w^{-1}-w^{-1},t)=y^{-1}w^{-1}-w^{-1}-t,
\end{equation}
which proves \eqref{eq:gen_eq}, so that $A(t)$ vanishes. This completes the proof.\qed

\section{The Abel Identities}\label{section:Abel}

In this section, we present a grammatical approach to the  Abel identites. To this end, we establish an expression of $D^n(ax^ry)$ in terms of
rooted trees on $[n]$. Recall that the set of such rooted
 trees is denoted by $R_n$.

For a rooted tree $T\in R_n$, we may construct a rooted tree $\bar{T}$
by coloring each child of the vertex $1$ by one of the colors $b,w_1,w_2,\ldots,w_r$. It should be noted that $1$ is not necessarily
the root of $T$.
Let $R^{(r)}_n$ denote the set of rooted trees on $[n]$
for which the children of $1$ are colored as described above.

We need the following grammatical labeling for a rooted tree $\bar{T}\in R^{(r)}_{n,k}$: First, represent $\bar{T}$ as a rooted tree $\hat{T}$ on $\{0,1,\ldots,n\}$ with root $0$, and represent an improper edge by double edges. Label the vertex $1$ by $ax^r$, label a black vertex by $x$ and label each of the remaining vertices by $w$. Moreover,
each edge in $\hat{T}$ is labeled by $y$. Thus the weight of $\bar{T}$ is given by
\begin{equation}\label{eq:wtrt}
w(\bar{T})=ax^{{\rm black}(\bar{T})+r}w^{n-1-{\rm black}(\bar{T})}y^{n+\rm imp(\bar{T})}.
\end{equation}

Using the same argument as in the proof of Theorem \ref{th:Dnab},
we are led to the following relation.

\begin{thm}
  For $n\geq1$ and $r\geq0$,
\begin{equation}\label{dn_weight_star}
  D^{n-1}(ax^ry)=\sum_{\bar{T}\in R^{(r)}_n}w(\bar{T}).
\end{equation}
\end{thm}

Analogous to Theorem \eqref{eq:dnaxr_close}, there is a connection between $D^n(ax^ry)$ and $Q_n(x,y)$.

\begin{thm}
  For $n\geq1$ and $r\geq0$,
\begin{equation}\label{eq:dna}
D^{n-1}(ax^ry)=ax^ry^{n}w^{n-1}Q_n(r+xw^{-1}-1,y).
\end{equation}
\end{thm}

In the notation of $Q_n(x,y)$, the relation \eqref{eq:Dumont1} of Dumont and Ramamonjisoa takes the form
\begin{equation}\label{eq:dnyw}
D^{n-1}(yw)=y^nw^{n}Q_n(0,y).
\end{equation}
In addition, Dumont \cite{Dumont} obtained grammatical expressions of $Q_n(1,y)$ and $Q_{n}(-1,y)$: For $n\geq1$,
\begin{align}
  D^{n}(w)&=y^{n}w^{n+1}Q_{n}(1,y),\label{eq:Dumont2}\\[6pt]
  D^{n}(y)&=y^{n+1}w^{n}Q_{n+1}(-1,y).\label{eq:Dumont3}
\end{align}
It can be checked that by setting $a=x=w$, the grammar $H$
 reduces to the grammar of Dumont and Ramamonjisoa. Meanwhile,  \eqref{eq:dnyw} can be deduced from \eqref{eq:dna} by setting $a=x=w$ and $r=0$ and  \eqref{eq:Dumont2} can be deduced from \eqref{eq:dnabr} by setting $a=x=w$ and $r=0$.

 We remark that \eqref{eq:Dumont3} can also be justified by
 a grammatical labeling of rooted trees in the set $R_{n,k}[\rm deg_T(1)=0]$
  of rooted trees in $R_{n,k}$  for which
   the vertex $1$ is a leaf.
 For a rooted tree $T$ in $R_{n,k}[\rm deg_T(1)=0]$, let  $\hat{T}$ denote the tree obtained from $T$ by adding a new root $0$ and representing each improper edge by double edges. Label each vertex except for $1$ by $x$ and label each edge in $\hat{T}$ by $y$. For a rooted tree in $R_{n,k}[\rm deg_T(1)=0]$,  it holds that
 \begin{equation}\label{eq:wTy}
   w(T)=y^{n+k}w^{n-1}.
 \end{equation}
 Since $1$ is not endowed with a label, this means that in Shor's construction, it is not allowed to add  new vertices as children of the vertex $1$. The argument for the proof of  Theorem \ref{grammar_w:Dnxy}
  implies that for $n\geq1$,
  \[
  D^n(y)=\sum_{T\in R_{n+1,k}[\deg_T(1)=0]}w(T).
  \]
  Utilizing the interpretation \eqref{zeng_2} of $Q_{n,k}(x)$, we see that
  for $x=-1$,
  \[
    Q_{n,k}(-1)=|R_{n+1,k}[\deg_T(1)=0]|.
  \]
  Thus it follows from \eqref{eq:wTy} that
  \[
    D^n(y)=y^{n+1}w^n\sum_{k=0}^{n-1}y^{k}Q_{n,k}(-1),
  \]
  which is the right hand side of \eqref{eq:Dumont3}.

The following   relations are needed in the
grammatical derivations of   Abel identities.

\begin{thm}
  For $n\geq1$,
  \begin{align}
    D^n(y)|_{y=w=1}&=n^n,\label{g:dny}\\[6pt]
    D^{n}(yw)|_{y=w=1}&=(n+1)^{n},\label{g:dnyw}\\[6pt]
    D^n(ax^r)|_{a=y=w=1}&=x^r(x+r)(x+r+n)^{n-1},\label{g:dnaxr}\\[6pt]
    D^n(ax^ry)|_{a=y=w=1}&=x^{r}(x+r+n)^{n}.\label{g:dnaxry}
  \end{align}
\end{thm}

\pf In the notation of $Q_n(x,y)$, the relation \eqref{Qnkx-sum} can be rewritten as
\begin{equation}\label{eq:sum_new}
Q_n(x,1)=(x+n)^{n-1}.
\end{equation}
Setting $y=w=1$ in \eqref{eq:Dumont3}, we obtain that
\[
D^n(y)|_{y=w=1}=Q_{n+1}(-1,1),
\]
which equals $n^n$ according to \eqref{eq:sum_new}. This proves  \eqref{g:dny}.
The rest of the relations in the theorem can be obtained from \eqref{eq:dnyw}, \eqref{eq:dnabr} and \eqref{eq:dna}, respectively.
This completes the proof. \qed

The classical Abel identity   states that for $n\geq1$,
\begin{equation}\label{eq:Abel_c}
(x+y+n)^n=\sum_{k=0}^n\binom{n}{k}x(x+k)^{k-1}(y+n-k)^{n-k}.
\end{equation}
Since $y$ has appeared as a variable in the grammar $H$, we shall use the
following form of \eqref{eq:Abel_c}: For $n\geq 1$,
  \begin{equation}\label{eq:Abel}
    (x_1+x_2+n)^n=\sum_{k=0}^n\binom{n}{k}x_1(x_1+k)^{k-1}(x_2+n-k)^{n-k}.
  \end{equation}

{\it\noindent Proof of \eqref{eq:Abel}.}  
Let
 \begin{align}
H'\colon &a_1\rightarrow a_1x_1y,\; a_2\rightarrow a_2x_2y,\; x_1\rightarrow x_1yw,\; x_2\rightarrow x_2yw,\notag\\[6pt]
&\qquad \qquad y\rightarrow y^3w,\; w\rightarrow yw^2,\label{g_d:g}
 \end{align}
and let $D$ denote the formal derivative  associated with the grammar  $H'$.
Viewing $a_1$ as $a$ and $x_1$ as $x$ and applying \eqref{g:dnaxr} with
 respect to $H$, we get
\begin{equation}\label{eq:dna1}
  D^n(a_1)|_{a_1=y=w=1}=x_1(x_1+n)^{n-1}.
\end{equation}
Similarly, invoking \eqref{g:dnaxry}, we obtain that
\begin{equation}\label{eq:dna2y}
  D^n(a_2y)|_{ a_2=y=w=1}=(x_2+n)^{n}.
\end{equation}
Moreover, since
\[
D(a_1a_2)=a_1a_2(x_1+x_2)y,\qquad D(x_1+x_2)=(x_1+x_2)yw,
\]
treating $a_1a_2$ as $a$ and $x_1+x_2$ as $x$, we may apply \eqref{g:dnaxry}  to deduce that
\begin{equation}\label{eq:dna12y}
  D^n(a_1a_2y)|_{a_1=a_2=y=w=1} = (x_1+x_2+n)^{n}.
\end{equation}
Finally, \eqref{eq:Abel} is follows from the Leibnitz formula
\begin{equation*}\label{eq:con_Abel1}
D^n(a_1a_2y)=\sum_{k=0}^n \binom{n}{k}D^k(a_1)D^{n-k}(a_2y)
\end{equation*}
along with the relations \eqref{eq:dna1}, \eqref{eq:dna2y} and \eqref{eq:dna12y}. This completes the proof. \qed

Riordan \cite{Riordan} obtained a class of generalizations of the classical Abel identity by considering the sum
\[
A_n(x,y;p,q)=\sum_{k=0}^n\binom{n}{k}(x+k)^{k+p}(y+n-k)^{n-k+q},
\]
where $n \geq 1$ and $p, q$ are integers.
He found closed forms for the cases when $(p,q)$ lies in
\begin{align*}
& \left\{(-3,0),(-2,0),(-1,0),(0,0),(1,0),\right.\\[6pt]
&\qquad \left.(2,0),(-1,-1),(-1,1),(-1,2),(1,1),(1,2),(2,2)\right\}.
\end{align*}
These identities for $A_n(x,y;p,q)$ are also called
the Abel identities or the Abel-type identities.
 The original Abel identity corresponds to the case $(p,q)=(-1,0)$. For $(p,q)=(-1,-1)$, Riordan obtained a closed form for $A_n(x,y;p,q)$, which is stated in the variables $x_1$ and $x_2$: For $n\geq1$,
  \begin{equation}\label{eq:Abel11}
    (x_1+x_2)(x_1+x_2+n)^{n-1}=\sum_{k=0}^n\binom{n}{k}x_1x_2(x_1+k)^{k-1}(x_2+n-k)^{k-1}.
  \end{equation}
{\it\noindent Proof of \eqref{eq:Abel11}.}
Let $H'$ denote the grammar given by \eqref{g_d:g}, and let $D$ denote the formal derivative associated with $H'$.
Using the same reasoning as for the proof of  \eqref{eq:dna1}, we see that
\begin{equation}\label{eq:dna2}
  D^n(a_2)|_{a_2=y=w=1}=x_2(x_2+n)^{n-1}.
\end{equation}
Analogous to \eqref{eq:dna12y}, we get the relation
\begin{equation}\label{eq:dna12}
  D^n(a_1a_2)|_{a_1=a_2=y=w=1}=(x_1+x_2)(x_1+x_2+n)^{n-1}.
\end{equation}
In view of \eqref{eq:dna1}, \eqref{eq:dna2} and \eqref{eq:dna12},
we are led to \eqref{eq:Abel11} by applying the
 Leibnitz formula
\[
D^n(a_1a_2)=\sum_{k=0}^n \binom{n}{k}D^k(a_1)D^{n-k}(a_2).
\]
This completes the proof. \qed

We next consider the case $(p,q)=(-2,0)$ of the Abel identity given by
Riordan: For $n\geq2$,
    \begin{align}
      &\sum_{k=0}^n\binom{n}{k}x_1(x_1+1)(x_1+k)^{k-2}(x_2+n-k)^{n-k}\notag\\[6pt]
      &\qquad =x_1^{-1}\left[(x_1+1)(x_1+x_2+n)^n-nx_1(x_1+x_2+n)^{n-1}\right].\label{eq:Moon3}
    \end{align}
{\it\noindent Proof of \eqref{eq:Moon3}.} Let $H'$ denote the grammar defined as above, and let $D$ denote the format derivative associated with $H'$. Set\
\[ s_1 = a_1y^{-1}+a_1x_1^{-1}w,\]
 so that \(D(s_1)=a_1x_1\). Analogous to \eqref{eq:dna1}, we find that for $n\geq1$,
\begin{equation}\label{eq:dna1x1}
  D^n(s_1)|_{a_1=y=w=1}=x_1(x_1+1)(x_1+n)^{n-2}.
\end{equation}
Since
\[
sa_2y= a_1a_2+a_1a_2x_1^{-1}yw
\]
and $x^{-1}w$ is a constant as shown in \eqref{eq:constantxw},
we get
\[ D^n(sa_2y)=D^n(a_1a_2)+x_1^{-1}wD^n(a_1a_2y).
\]
By the Leibnitz formula
\[
 D^n(sa_2y)=\sum_{k=0}^n\binom{n}{k}D^k(s)D^{n-k}(a_2y),
\]
we find that
\[
 \sum_{k=0}^n\binom{n}{k}D^k(s)D^{n-k}(a_2y)=D^n(a_1a_2)+x^{-1}wD^n(a_1a_2y),
\]
which yields \eqref{eq:Moon3} by applying  \eqref{eq:dna2y},  \eqref{eq:dna12y}, \eqref{eq:dna12} and \eqref{eq:dna1x1}. This completes the proof.  \qed

Next we provide a closed formula
 for the case $(p,q)=(-2,-2)$.

\begin{thm}\label{th:Abel22}
  For $n\geq1$,
  \begin{align}
    &\sum_{k=0}^n\binom{n}{k}x_1x_2(x_1+1)(x_2+1)(x_1+k)^{k-2}(x_2+n-k)^{n-k-2}\notag\\[6pt]
    &\qquad=((x_1+x_2)^3-3n(x_1+x_2)-2n)(x_1+x_2+n)^{n-3} \notag\\[3pt]
    &\qquad\qquad+\frac{(x_1+x_2)^2}{x_1x_2}(x_1+x_2+1)(x_1+x_2+n)^{n-2}.\label{eq:Abel22}
  \end{align}
\end{thm}

We need the following lemma.

\begin{lem}
Let $D$ denote the formal derivative associated with the grammar $H$. For $n\geq2$,
  \begin{equation}\label{g:dnax1}
    D^n(ay^{-1})|_{a=y=w=1}=x(x+1)(x+n)^{n-2}-(x+n)^{n-1}.
  \end{equation}
For $n\geq3$,
  \begin{equation}\label{g:dnax2}
    D^n(ay^{-2})|_{a=y=w=1}=(x^3-3nx-2n)(x+n)^{n-3}.
  \end{equation}
\end{lem}
\pf   Since $D(ay^{-1})=ax-ayw$ and $D(x^{-1}w)=0$ as given in \eqref{eq:constantxw}, we have
  \[
    D^n(ay^{-1})=D^{n-1}(ax)-x^{-1}wD^{n-1}(axy).
  \]
  Applying \eqref{g:dnaxr} and \eqref{g:dnaxry} with $r=1$, we get
  \[
  D^n(ay^{-1})|_{a=y=w=1}=x(x+1)(x+n)^{n-2}-(x+n)^{n-1}.
  \]
Since
  \[
    D(ay^{-2})=axy^{-1}-2aw
  \]
  and
  \[
    D(axy^{-1})=ax^2(1+x^{-1}w)-axyw,
  \]
 we see that  for $n\geq3$,
  \begin{align*}
    D^n(ay^{-2})&=D^{n-1}(axy^{-1}-2aw)\\[6pt]
    &=D^{n-2}\left(ax^2(1+x^{-1}w)-axyw\right)-2x^{-1}wD^{n-1}(ax)\\[6pt]
    &=(1+x^{-1}w)D^{n-2}(ax^2)-x^{-1}wD^{n-2}(ax^2y)-2x^{-1}wD^{n-1}(ax).
  \end{align*}
 In light of \eqref{g:dnaxr} and \eqref{g:dnaxry}, we find that
  \[
      D^n(ay^{-2})|_{a=y=w=1}=x(x+1)(x+2)(x+n)^{n-3}-x(x+n)^{n-2}
      -2(x+1)(x+n)^{n-2},
  \]
  which implies \eqref{g:dnax2}. This completes the proof. \qed

{\noindent \it Proof of Theorem \ref{th:Abel22}.  }  Assume that
 $H'$  the grammar given in \eqref{g_d:g} and $D$ is the
 formal derivative with respect to $H'$. Let
\[
s_1=a_1y^{-1}+a_1x_1^{-1}w
\]
and
\[
s_2=a_2y^{-1}+a_2x_2^{-1}w.
\]
Clearly,
\(
D(s_1)=a_1x_1
\)
and
\(
D(s_2)=a_2x_2.
\)
It follows from \eqref{g:dnaxr} that
\begin{equation}\label{eq:s2}
  D^n(s_2)|_{a_2=y=w=1}=x_2(1+x_2)(x_2+n)^{n-2}.
\end{equation}
By the same argument  as that for the proof of \eqref{eq:dna12y},
 we deduce from  \eqref{g:dnax1} and \eqref{g:dnax2} that
\begin{align*}
    D^n(a_1a_2y^{-1})|_{a_1=a_2=y=w=1}
  &=(x_1+x_2)(x_1+x_2+1)(x_1+x_2+n)^{n-2}\\[6pt]
  &\qquad -(x_1+x_2+n)^{n-1}
  \end{align*}
  and
  \[
  D^n(a_1a_2y^{-2})|_{a_1=a_2=y=w=1}
  =((x_1+x_2)^3-3n(x_1+x_2)-2n)(x_1+x_2+n)^{n-3}.
\]
Since $D(x_1w^{-1})=D(x_2w^{-1})=0$, we get
\begin{align}
  D^n(s_1s_2)&=D^n(a_1a_2(y^{-1}+x_1^{-1}w)(y^{-1}+x_2^{-1}w))\notag\\[6pt]
  &=D^n\left(a_1a_2y^{-2}+a_1a_2y^{-1}(x_1^{-1}+x_2^{-1})w+a_1a_2(x_1x_2)^{-1}w^2\right)\notag\\[6pt]
  &=D^n(a_1a_2y^{-2})+(x_1^{-1}  +x_2^{-1})wD^n(a_1a_2y^{-1})\notag\\[6pt]
  &\qquad +(x_1x_2)^{-1}w^2D^n(a_1a_2).\label{eq:s12}
\end{align}
By  the Leibnitz formula
\[
    D^n(s_1s_2)=\sum_{k=0}^n\binom{n}{k}D^k(s_1)D^{n-k}(s_2),
\]
we obtain \eqref{eq:Abel22} by using \eqref{eq:dna1x1}, \eqref{eq:s2} and \eqref{eq:s12}. This completes the proof.
\qed

We conclude this paper with a one-line grammatical
explanation of  the Lacasse identity.

\begin{thm}
  For $n\geq 1$,
  \begin{equation}\label{eq:ML}
  n^{n+1}=\sum_{j=1}^n\sum_{k=0}^{n-j}\binom{n}{j}\binom{n-j}{k}j^jk^k(n-j-k)^{n-j-k}.
  \end{equation}
\end{thm}
\pf Because of the relation
\[
k\binom{n}{k}=n\binom{n-1}{k-1},
\]
\eqref{eq:ML} can be rewritten as
\begin{equation}\label{eq:ML_eq}
n^n=\sum_{j=1}^n\sum_{k=0}^{n-j}\binom{n-1}{j-1}\binom{n-j}{k}j^{j-1}k^k(n-k)^{n-k}.
\end{equation}
Since $D(y)=y^3w$, we have
\begin{equation}\label{eq:cony}
D^{n}(y)=D^{n-1}(y^3w)=\sum_{i+j+k=n-1}\binom{n-1}{i,j,k}D^{i}(y)D^{j}(yw)D^{k}(y).
\end{equation}
Invoking \eqref{g:dny} and \eqref{g:dnyw} and setting $y=w=1$, we see that \eqref{eq:cony} can be rewritten in the form
 of \eqref{eq:ML_eq}. This completes the proof. \qed

\vspace{.2cm} \noindent{\bf Acknowledgment.}
This work was supported by the National
Science Foundation of China.


\begin{thebibliography}{99}

\bibitem{Berndt} B. C. Berndt, R. J. Evans and B. M. Wilson, Chapter 3 of Ramanujan's second notebooks, Adv. Math. 49 (1983) 123--169.

\bibitem{Chen} W.Y.C. Chen, Context-free grammars, differential
operators and formal power series, Theoret. Comput. Sci. 117 (1993) 113--129.

\bibitem{Fu} W.Y.C.~Chen and A.M.~Fu, Context-free grammars, permutations and increasing trees, Adv. in Appl. Math. 82 (2017) 58--82.

\bibitem{CG} W.Y.C Chen and V.J.W Guo, Bijections behind the Ramanujan polynomials, Adv. in Appl. Math. 27 (2001) 336--356.

\bibitem{CPY} W.Y.C. Chen, J.F.F. Peng, H.R.L. Yang, Decomposition of triply rooted trees, Electron. J. Combin. 20(2) (2013) \#P10.

\bibitem{dumont} D. Dumont, Grammaires de William Chen et d\'{e}rivations dans les arbres et
arborescences, S\'{e}m. Lothar. Combin. 37 (1996) B37a, 21 pp.

\bibitem{Dumont} D. Dumont and A. Ramamonjisoa, Grammaire de Ramanujan et Arbres de Cayley,
Electron. J. Combin. 3 (1996) R17.

\bibitem{Gessel} I.M. Gessel, Lagrange inversion, J. Combin. Theory Ser. A 144 (2016) 212--249.

\bibitem{Guo}V.J.W. Guo, A bijective proof of the Shor recurrence, European J. Combin. 70 (2018) 92--98.


\bibitem{Lacasse} A. Lacasse, Bornes PAC-Bayes et algorithmes d¡¯apprentissage, Ph.D. Thesis, Universite Laval, Quebec, 2010.

\bibitem{MMY}S.M. Ma, J. Ma and Y.N. Yeh. The ascent-plateau statistics on Stirling permutations,  arXiv:1801.08056, 2018.

\bibitem{MY}S.M. Ma and Y.N. Yeh, Eulerian polynomials, perfect matchings and Stirling permutations of the second kind, Electron. J. Combin. 24  (4) (2017) \#P27

\bibitem{Prodinger} H. Prodinger,  An identity conjectured by Lacasse via the tree function,  Electron. J. Combin. 20 (3) (2013) \#P7.

\bibitem{Ramanujan}S. Ramanujan, Notebooks, Vol. 1, Tata Institue of Fundanmental Research, Bombay, 1957, 35--36

\bibitem{Riordan}J. Riordan, Combinatorial Identities, New York: Wiley, 1968.

\bibitem{Rota}G.-C. Rota, Finite Operator Calculus, Academic Press, New York, 1975.

\bibitem{Shor} P.W. Shor, A new proof of Cayley's formula for counting labelled trees. J. Combin. Theory Ser. A 71 (1995) 154--158.

\bibitem{Sun} Y. Sun, A simple proof of an identity of Lacasse,  Electron. J. Combin. 20(2) (2013) \#P11.

\bibitem{Younsi} M. Younsi, Proof of a combinatorial conjecture coming from the PAC-Bayesian machine learning theory, arXiv:1209.0824, 2012.

\bibitem{Zeng} J. Zeng, A Ramanujan sequence that refines the Cayley formula for trees,
Ramanujan J. 3 (1999) 45--54.
\end{thebibliography}
\end{document}